\documentclass[12 pt]{amsart}

\usepackage{amssymb,amsfonts,graphicx}
\usepackage{yhmath}
\usepackage[all,arc]{xy}
\usepackage{enumerate}
\usepackage{mathrsfs}
\usepackage{enumitem} 
\usepackage{url}

\newtheorem{theorem}{Theorem}[section]
\newtheorem{corollary}[theorem]{Corollary}

\newtheorem{lemma}[theorem]{Lemma}

\theoremstyle{definition}
\newtheorem{definition}[theorem]{Definition}

\newtheorem{example}[theorem]{Example}
\newtheorem{counterexample}[theorem]{Counterexample}

\newtheorem{notations}[theorem]{Notations}

\theoremstyle{plain}

\DeclareMathAlphabet{\mathpzc}{OT1}{pzc}{m}{it}

\newcommand{\N}{\mathbb{N}}
\newcommand{\R}{\mathbb{R}}
\newcommand{\C}{\mathbb{C}}

\theoremstyle{definition}

\newtheorem{remark}[theorem]{Remark}

\newcommand{\diff}{\mathrm{d}}

\providecommand{\C}{\operatorname{\mathcal{C}^*}}

\DeclareMathAlphabet{\mathpzc}{OT1}{pzc}{m}{it}

\theoremstyle{remark}

\providecommand{\C}{\operatorname{Cap}}

\DeclareMathAlphabet{\mathpzc}{OT1}{pzc}{m}{it}


\numberwithin{equation}{section}

\title[Subharmonic and Plurisubharmonic Functions]{Families of (Pluri)subharmonic Functions and Separate Subharmonicity}

\author[M. Kalantar]{Mansour Kalantar}

\begin{document}
	
	\begin{abstract}
		We study the singular sets of pointwise bounded-above envelopes of families of functions. For families of subharmonic functions on open subsets of $\mathbb C$, of plurisubharmonic functions on open subsets of $\mathbb C^N$, and, more generally, of functions subharmonic in each complex variable separately, we prove that the singular set of a pointwise bounded-above envelope is relatively closed and nowhere dense. In contrast, for $N\ge3$ we construct a sequence of subharmonic functions on $\mathbb R^N$ whose upper envelope is finite everywhere but nowhere locally bounded above. We also characterize the singular sets that arise for admissible classes of functions, both with and without a continuity assumption, and give applications to separately subharmonic functions.

	\end{abstract}
	\maketitle

	\section{Introduction}
	Let $\Omega$ be a domain in $\R^N$ (resp., $\C^N$), $N\ge2$,  and let
	$\{{u_\alpha}\}_{\alpha\in J}$ be a nonempty family of
	subharmonic (resp., plurisubharmonic) functions on $\Omega$. Consider its upper envelope
	$
	u(z):=\sup_{\alpha\in J}u_\alpha(z)
	$
	and its upper semicontinuous regularization
	$$
	u^*(z):=\limsup_{\substack{\zeta\to z\\(\zeta\in\Omega)}}u(\zeta).
	$$
	A classical theorem in potential (resp., pluripotential) theory asserts that, if the
	family is locally uniformly bounded above, then $u^*$ is
	subharmonic (resp., plurisubharmonic) and  the exceptional set
	$$
	\{z\in\Omega: u(z)<u^*(z)\}
	$$
	is polar (resp., pluripolar). 
	
	It is therefore natural to ask under what conditions the upper
	envelope of a family of subharmonic or plurisubharmonic functions is locally bounded
	above. Pointwise finiteness,
	$$
	u(z)<+\infty\qquad(z\in\Omega),
	$$
	does not by itself imply local boundedness above
	(see Theorem \ref{S_Can_be_R}). This leads us to study the set of points at which the upper envelope fails to be locally bounded above, namely the singular set. Also, the case of families of subharmonic functions has applications to separate subharmonicity.
	
	Kosiński, Martel, and Ransford \cite{KosinskiMartelRansford} proved a uniform boundedness theorem for families $\mathcal U$ of continuous plurisubharmonic functions on $D\times G$, where $D\subset\mathbb C^N$ and $G\subset\mathbb C^M$ are domains and $N,M\ge 1$. According to their result, if $\mathcal U$ is uniformly bounded above on $K\times\{w\}$ for every compact $K$ of $ D$ and every $w\in G$, and on $\{z\}\times L$ for every $z\in D$ and every compact $L$ of $ G$, then $\mathcal U$ is uniformly bounded above on every compact subset of $D\times G$.

	In this paper, we prove that a pointwise bounded-above upper envelope of subharmonic functions on an open subset of \(\mathbb C\) is locally bounded above outside a relatively closed nowhere dense set. The same conclusion holds for families of plurisubharmonic functions defined on an open subset $\C^N$, $N\ge1$, and, more generally, for families of functions that are subharmonic in each complex variable separately. For \(N\geq3\), we construct a sequence of subharmonic functions on \(\mathbb R^N\) whose upper envelope is finite everywhere but nowhere locally bounded above. We then characterize the singular sets that arise for admissible classes of functions, both with and without a continuity assumption. We also give a geometric sufficient condition for a closed subset of \(\mathbb R^N\) to arise as such a singular set for a family of subharmonic functions. We also apply our results to separately subharmonic functions on product domains.

	\section{Preliminaries}
	
		\begin{notations}
		We write $\N=\{1,2,\dots\}$ for the set of natural numbers.  For a set $A$ in the ambient topology $X$, we denote its closure and interior in $X$ by $\overline A$ and $\operatorname{int}_X(A)$, respectively. We say $A$ is relatively compact in $X$ if the closure of $A$ in $X$ is compact. Euclidean balls are written $B(a,r)$ and the ambient topology is clear from the context. 
		
	\end{notations}
	
	Let $\Omega$ be an open subset of $\R^N$, with $N\ge2$. We say an upper semicontinuous function
	$
	u:\Omega\to[-\infty,+\infty)
	$
	is \textit{subharmonic}, and write $u\in\mathcal{SH}(\Omega),$ if $u\not\equiv-\infty$ on each component of $\Omega$, and for every open ball $B(x_0,r)$ that is relatively compact in $\Omega$,
	\[
	u(x_0)\le \frac{1}{|B(x_0,r)|}
	\int_{B(x_0,r)} u(\xi)\diff\mu(\xi),
	\]
	where $\mu$ is Lebesgue measure in $\R^N$ and $|B(x_0,r)|$ denotes the Lebesgue measure of the ball.

	Now, let $\Omega$ be an open subset of $\C^N$, $N\geq 1$.
	We say that a function $u:\Omega\to[-\infty,+\infty)$ is
	\textit{plurisubharmonic}, and write $u\in\mathcal{PSH}(\Omega)$, if $u$ is not identically equal to $-\infty$ on any component of $\Omega$, is upper semicontinuous and for every complex affine line $L$ intersecting $\Omega$, the restriction $u\big|_{\Omega\cap L}$ is subharmonic or identically $-\infty$ on every component of the set $\Omega\cap L$ (see Klimek \cite{Klimek}).

	In this paper, $\Omega$ denotes an open subset of $\R^{N}$,
	$N\ge2$, in the subharmonic setting, or of $\C^{N}$,
	$N\ge1$, in the plurisubharmonic setting.

	\begin{theorem}[El Marzguioui--Wiegerinck]
		Let $U\subset\C^N$ be plurifinely open and let $a\in U$.
		Then there exist $\kappa>0$ and a plurifine neighborhood
		$V\subset U$ of $a$ such that, for every $v\in V$ and every complex
		line $L$ passing through $v$, the plurifine component of $U\cap L$
		containing $v$ contains a circle centered at $v$ whose radius is at	least $\kappa$.
	\end{theorem}
	See \cite[Proposition 4.1]{ElMarzguiouiWiegerinck2009}.
	We will use the following version of this result.
	
	\begin{corollary}\label{El_Marz_Wieg Theorm}
		Let $\Omega\subset\mathbb C^N$ be open, and let $U\subset\Omega$
		be a plurifine neighborhood of a point $a\in\Omega$. Then there exists $r>0$
		such that  the ball $B(a,r)$ is relatively compact in $\Omega$ and  for every
		$v\in\mathcal{PSH}(\Omega)$ and every $M\in\mathbb R$,
		\[
		v\le M \text{ on }U
		\quad\Longrightarrow\quad
		v\le M \text{ on }B(a,r).
		\]
	\end{corollary}
	
	\begin{proof}
		Choose  a Euclidean ball $B(a,\rho)$ centered at  $a$ and relatively compact in $\Omega$, and a plurifinely open
		neighborhood $W$ of $a$ contained in $U\cap B(a,\rho)$.
		By El Marzguioui--Wiegerinck's theorem, there exists
		$\kappa>0$ such that, for every complex line $L$ through $a$,
		there is a radius $r_L\ge\kappa$ satisfying
		$\partial_L(B(a,r_L)\cap L)\subset W$, where $\partial_L$
		denotes the boundary relative to $L$. Since the circle lies in \(B(a,\rho)\), its radius \(r_L\) is less than \(\rho\); hence its enclosed disc is relatively compact in \(\Omega\cap L\).
		If $v\le M$ on $U$, the maximum principle gives $v\le M$
		throughout each such disc. These discs cover $B(a,\kappa)$,
		so $v\le M$ there. Taking $r=\kappa/2$ completes the proof.
	\end{proof}

	Our reference for what follows in the real case is Armitage and Gardiner \cite[Chapter 7]{Armitage Gardiner book} and in the complex case are Klimek \cite{Klimek} and  El Kadiri--Fuglede \cite{ElKadiri_Fuglede}.  The \textit{fine topology} (resp., \textit{plurifine topology})on $\C\cong \R^2$ (resp., $\C^N$, $N\ge1$) is the coarsest topology on $\C$ (resp., $\C^N$) that makes all subharmonic (resp., plurisubharmonic) functions continuous in the extended sense. As a consequence, each subharmonic (resp., plurisubharmonic) function defined on an open set $U$ of $\C$ (resp., of $\C^N$) is finely (resp., plurifinely) continuous on $U$.
	We use the adjective \textit{fine} (resp., plurifine) and the adverb \textit{finely} (resp., plurifinely) to refer to the fine (resp., plurifine) topology. When needed, a subscript “f” (resp., $\mathcal{F}$) is used to indicate that a notation is taken with respect to the fine (resp., plurifine) topology.
	
	A set $E\subset\C$ is said to be \textit{thin} at a point $z\in\C$ if $z$ is not a fine limit point of $E$. In other words, there is a fine neighborhood $V$ of $z$ which does not intersect $E\setminus\{z\}$. A set \(V\subset\mathbb{R}^N\) containing $x$ is a fine neighborhood of $x$ if its complement $\R^N\setminus V$ is thin at $x$ (\cite[Theorem 1.2.4, p. 13]{ElKadiri_Fuglede}). 
	
	 If a set $E\subset\C$ is thin at $z$, then there are arbitrarily small circles centered at $z$ which do not intersect $E$ (see Armitage and Gardiner \cite[Theorem 7.3.9]{Armitage Gardiner book}). This result does not hold in $\R^N$ for $N\ge3$. 
	One application of this result is the following theorem that we use frequently. For a stronger version, see Fuglede, Groot, and Wiegerinck \cite[Lemma~2.4]{FugledeGrootWiegerinck2019}.
	
	\begin{lemma}[Fine-neighborhood circle lemma]\label{Fine-neigh_circl}
		Let $V$ be a fine neighborhood of $z_0\in\C$, and let
		$\mathcal{O}$ be a Euclidean neighborhood of $z_0$. Then there exists
		$r>0$ such that $D(z_0,r)$ is relatively compact in $\mathcal{O}$ and	$\partial D(z_0,r)\subset V$.
		
	\end{lemma}
	\begin{proof}
		Since $V$ is a fine neighborhood of $z_0$, its complement
		$\C\setminus V$ is thin at $z_0$.
		By the above result, there are
		arbitrarily small circles centered at $z_0$ that do not
		intersect $\C\setminus V$. Choose $\rho>0$ such that
		$D(z_0,\rho)$ is relatively compact in $\mathcal{O}$. We may therefore choose $0<r<\rho$ such that
		\[
		\partial D(z_0,r)\cap (\C\setminus V)=\varnothing.
		\]
		Thus $\partial D(z_0,r)\subset V$ and
		$D(z_0,r)$ is relatively compact in $\mathcal{O}$, as required.
		
	\end{proof}
	
	\begin{theorem}\label{RungeType} 
		Let $F$ be a compact subset of an open set $\Omega$ in  $ \R^{N}$ such that every bounded component of $\R^{N} \setminus F$ contains a point of $\R^{N}\setminus\Omega$. If $h$ is harmonic on an open set containing $F$ and if $\varepsilon>0$, then there exists a function $H$ harmonic on $\Omega$ such that $|h-H|<\varepsilon$ on $F$.
	\end{theorem} 
	See Armitage and Gardiner\cite[Theorem 2.6.4.]{Armitage Gardiner book}.

	\section*{Main Results}

	Throughout this paper, $\Omega$ denotes an open subset of
	$\C\cong\R^2$, $\C^N$
	with $N\geq 1$, or $\R^N$ with $N\geq 2$.
	We consider a family $\{u_\alpha\}_{\alpha\in J}$ of functions
	from $\Omega$ to $[-\infty,+\infty)$ and set
	\[
	u:=\sup_{\alpha\in J}u_\alpha.
	\] 
	\begin{definition}\label{regular_set_singular_set}
		The \textit{regular set }$R_u$ of $u$ consists of all points
		$\zeta\in\Omega$ at which $u$ is locally bounded above.
		The \textit{singular set} of $u$ is $S_u:=\Omega\setminus R_u$.
	\end{definition}
	
	By definition, $R_u$ is open and $S_u$ is relatively closed in $\Omega$. Saying that $R_u$ is dense in $\Omega$ is equivalent to saying that $S_u$ is nowhere dense in $\Omega$. We shall drop the subscript $u$ if there is no confusion.

	\section{Families of Plurisubharmonic Functions}
	
	\begin{theorem}
		Let $\Omega$ be an open subset of $\C^N$, $N\ge1$, and let
		$\{u_\alpha\}_{\alpha\in J}$ be a nonempty family of
		plurisubharmonic functions on $\Omega$. Suppose that
		\begin{equation}\label{hypothesis_Theo_3_1}
			u(z):=\sup_{\alpha\in J}u_\alpha(z)<+\infty
		\end{equation}
			for every $	z\in\Omega$.
		Then the singular set $S$ of $u$ is a relatively closed nowhere dense subset of $\Omega$.
		Moreover, the upper semicontinuous regularization $u^*$ of
		$u$ is plurisubharmonic on
		$\Omega\setminus S$ and the negligible set
		\begin{equation}\label{negligible_set_is_pluripolar}
		\{z\in
		\Omega\setminus S: u(z)<u^*(z)\}
		\end{equation} is pluripolar in $\C^N$.
		
	\end{theorem}
	
	\begin{proof}
		Denote the plurifine topology by $\mathcal{F}$. For $m\in\N$, let 
		$	E_m$ be the set of all $z\in\Omega$ such that $u(z)\le m$. Since plurisubharmonic functions are  plurifinely continuous, $u$ is plurifinely lower semicontinuous, and each $E_m$ is plurifinely closed.  Using \eqref{hypothesis_Theo_3_1}, one checks that
		\[	\Omega=\bigcup_{m=1}^\infty E_m.
		\] Since the plurifine topology in $\C^N$ is a Baire topology (\cite[p. 323]{ElKadiri_Fuglede}), it follows that the set 
		\[
		G:=\bigcup_{m=1}^{\infty}\operatorname{int}_{\mathcal{F}}(E_m)
		\]  is plurifinely dense in $\Omega$, and therefore also Euclidean	dense.
		
	To see that $u$ is Euclidean locally bounded on $G$, take a point $a\in G$. Choose $m$ such that
	$a\in U:=\operatorname{int}_{\mathcal F}(E_m)$. 
	The set $U$ is a plurifine neighborhood of $a$, and $u\le m$ on $U$. This implies
	\begin{equation}\label{u_less_than_m}
	u_\alpha(z)	\le m
	\end{equation}
	for all  $z\in U$ and $\alpha\in J$. Applying 
	Corollary \ref*{El_Marz_Wieg Theorm}, we obtain a radius $r>0$ independent of $\alpha$ such that the ball $B(a,r)$ is relatively compact in $\Omega$ and on which inequality \eqref{u_less_than_m} holds for all $\alpha\in J$.  Taking the supremum over $\alpha$, we obtain 
	\[ u(z)\le m\] for all $z\in B(a,r)$.
	Therefore $u$ is locally bounded at each point of $G$.
	
	 Consequently,  the
		regular set $R$ of $u$ contains $G$. Since according to the first paragraph $G$ is Euclidean dense in $\Omega$, we obtain that $R$ is Euclidean dense in $\Omega$, as well.   Thus the singular set
		$S:=\Omega\setminus R$, which  is relatively closed in $\Omega$, is nowhere
		dense and $u$ is locally bounded above on $\Omega\setminus S$.
		
		Finally, since $u$ is locally bounded above on $\Omega\setminus S$, classical results show that its upper semicontinuous regularization is plurisubharmonic there \cite[Theorem~2.9.14-(iv), p. 69]{Klimek}, and its negligible set, defined in \eqref{negligible_set_is_pluripolar}, is pluripolar in $\C^N$ (\cite[p. 166 and Theorem 4.7.6, p. 171]{Klimek}).

	\end{proof}

	The next result follows from the last theorem and the fact that the classes of subharmonic functions and plurisubharmonic functions are identical in $\C$.  We identify $\C$ with $\R^2$ via the map
	$x+iy \mapsto (x,y)$.

	\begin{corollary} \label{PlanarUniformBT}
		Let $\Omega\subset\C$ be open, and let
		$\{u_\alpha\}_{\alpha\in J}$ be a family of subharmonic functions on
		$\Omega$. Suppose that
		\[
		u(z):=\sup_{\alpha\in J}u_\alpha(z)<+\infty
		\qquad (z\in\Omega).
		\]
		Then the singular set $S$ of $u$ is a relatively closed nowhere dense subset of $\Omega$.
	\end{corollary}

	
	\section{Families of Subharmonic Functions in $\R^N$ for $N\ge3$}

	The analogue of the previous corollary fails in $\R^N$ for $N\ge3$, as the following counterexample shows.
	
	\begin{theorem}\label{S_Can_be_R}
		Let $N\ge3$. There exists a sequence $(u_n)_{n\in\N}$ of nonnegative subharmonic functions on $\R^{N}$ such
		that	\[
		u(z):=\sup_{n\in\N}u_n(z)<+\infty
		\] for every $	x\in\R^N$,
		but $u$ is not locally bounded above at any point of $\R^N$.
	\end{theorem}
	
	\begin{proof}

\leavevmode\par
\medskip
\noindent\textbf{Step 1}\par

		Let
		\[
		D=
		\bigcup_{j=1}^{N}
		\left\{
		x=(x_1,\ldots,x_N)\in\R^{N}:
		x_i\in\mathbb{Q}\text{ for every }i\ne j
		\right\}.
		\]
		The set $D$ is dense because it contains $\mathbb{Q}^{N}$. It is also path connected. Indeed, let $x\in D$. After relabeling the coordinates, assume that \(x_2,\ldots,x_N\) are rational. The segment
		\[
		\bigl\{(t x_1,x_2,\ldots,x_N):0\le t\le1\bigr\}
		\]
		joins \(x\) to \((0,x_2,\ldots,x_N)\) and lies in \(D\), because all coordinates except the first are rational.
		Once the first coordinate is $0$, all coordinates are rational and we may vary $x_2$ in $\R$ while remaining in  $D$. We then use the segment
		\[
		\bigl\{(0,t x_2,x_3,\ldots,x_N):0\le t\le1\bigr\}
		\]
		to reach \((0,0,x_3,\ldots,x_N)\). This segment also lies in \(D\), because all coordinates except the second are rational.  Continuing in this way, we reach the origin after \(N\) segments, without leaving \(D\). Thus every point of \(D\) can be joined to the origin by a continuous path in \(D\). Thus \(D\) is path connected and so connected.

		The set $D$ is also polar. For \(j=1,\ldots,N\) and \(\alpha=(\alpha_1,\ldots,\alpha_{N-1})\in\mathbb Q^{N-1}\), let
		\[
		L_{j,\alpha}
		:=
		\{(\alpha_1,\ldots,\alpha_{j-1},t,\alpha_j,\ldots,\alpha_{N-1}):t\in\R\}.
		\]
		Then
		\[
		D=\bigcup_{j=1}^{N}\ \bigcup_{\alpha\in\mathbb Q^{N-1}}L_{j,\alpha}.
		\]
		Every affine subspace of dimension at most \(N-2\) is polar in \(\R^N\), and \(L_{j,\alpha}\) has dimension 1. Thus each line \(L_{j,\alpha}\) is polar in \(\R^N\) when \(N\ge3\).
		Since \(\mathbb Q^{N-1}\) is countable and countable unions of polar sets are polar, \(D\) is polar.
		
		It follows from a standard theorem (see   \cite[Theorem 5.1.3(i), p. 124]{Armitage Gardiner book}) that there exists a nonnegative superharmonic function $p$ on $\R^{N}$ such that $p\not\equiv+\infty$ and $p=+\infty$ on $D$.
		
		\medskip
		\noindent
		\textbf{Step 2}
		
		Let $|\cdot|$ designate the Euclidean norm in $\R^N$. For $n\in\N$, define
		\[
		K_n=
		\left\{x\in\R^{N}:
		|x|\le n,\ p(x)\le n
		\right\}.
		\]
		Take $|\zeta_0|\le 1$ such that $
		m:=\min_{|x|\le1}p(x)=p(\zeta_0).$
		Replacing $p$ by $p/\max\{1,m\}$, we may assume that
		$m\le1$. This normalization preserves nonnegativity,
		superharmonicity, and the set where $p=+\infty$.
		Consequently $K_1$ is not empty because it contains $\zeta_0$. Since the sequence of sets $(K_n)_{n\in\N}$ is increasing, $K_n$ is nonempty 
		for every $n\in\N$.
		Also, since $p$ is lower semicontinuous, $\{p\le n\}$ is closed,
		so $K_n$ is compact. One checks that the increasing sets $K_n$ satisfy 
		\[
		\bigcup_{n\ge1}K_n=\{x:p(x)<+\infty\}.
		\] 
		Moreover, since $p=+\infty$ on $D$, we obtain that $K_n$ does not intersect $D$, and so
		\begin{equation}\label{K_n_is_inside_comple_D}
				 K_n\subset\R^N\setminus D.
		\end{equation}
	 Since $D$ is dense, each
		$K_n$ has empty interior:
		\begin{equation}\label{int_K_n=_empty}
			\operatorname{int}(K_n)=\emptyset.
		\end{equation}
		
		We shall repeatedly use the following observation:
		if a compact set $E$ is disjoint from $D$, then
		$\R^{N}\setminus E$ is connected.
		Indeed, this complement contains $D$, which is connected
		and dense in the complement. A space containing a dense
		connected subset is connected. This, applied to $K_n$, gives $\R^N\setminus K_n$ is connected. 
		
		\medskip
		\noindent
		\textbf{Step 3}
		
		The balls with rational centers and positive rational
		radii form a countable basis of $\R^{N}$.
		Enumerate them as
		\[
		B_1,B_2,B_3,\ldots.
		\]
		Define a sequence $(V_n)_{n\in\N}$ by concatenating the finite lists
		\[
		(B_1),\quad
		(B_1,B_2),\quad
		(B_1,B_2,B_3),\quad\ldots.
		\]
		Thus every fixed basic ball occurs infinitely often,
		at arbitrarily large indices. For each $n$, the set $V_n\setminus K_n$ is open and nonempty because $K_n$ has empty interior by \eqref{int_K_n=_empty}. Since $p$ is finite
		almost everywhere, choose
		$b_n\in V_n\setminus K_n$ with $p(b_n)<+\infty$.
		In particular,
		\begin{equation}\label{b_n_notin_D}
			b_n\notin D.
		\end{equation}
		\medskip
		\noindent
		\textbf{Step 4}
		
		Set
		\[
		E_n=K_n\cup\{b_n\}.
		\]
		This compact set is disjoint from $D$ by (\ref{K_n_is_inside_comple_D}) and (\ref{b_n_notin_D}), so its complement
		is connected by the remark at the end of Step 2. Since $b_n\notin K_n$, choose disjoint open neighborhoods
		of $K_n$ and $\{b_n\}$. Define a harmonic function that equals $0$ on a small neighborhood of $K_n$ and equals $p(b_n)+n+2 $ on a small neighborhood of $\{b_n\}$.
		  Theorem \ref{RungeType} provides a harmonic function $h_n$ on $\R^N$ satisfying
		\begin{equation}\label{properties_of_h_n}
			|h_n(x)|<\frac12
		 \textrm{ for } x\in K_n\quad\textrm{ and}\quad
		h_n(b_n)>p(b_n)+n+\frac32.
		\end{equation}

		\medskip
		\noindent
		\textbf{Step 5}
		
		Define
		\[
		u_n(x)=\max\{h_n(x)-p(x)-1,0\}.
		\]
		Since $-p$ is subharmonic and $h_n$ is harmonic,
		$h_n-p-1$ is subharmonic, and so is $u_n$. 
		Thus every $u_n$ is nonnegative and
		subharmonic on $\R^{N}$.
		For $x\in K_n$, we have by \eqref{properties_of_h_n} and in view of the nonnegativity of $p$,
		\[
		h_n(x)-p(x)-1
		<\frac12-0-1<0,
		\]
		and hence 
		\begin{equation}\label{u_n=0_on_K_n}
			u_n=0 \textrm{  on } K_n
		\end{equation}
		and
		\begin{equation}\label{u_n>n_at_b_n}
		u_n(b_n)
		=h_n(b_n)-p(b_n)-1
		>n+\frac12>n.
			\end{equation}
		
		\medskip
		\noindent
		\textbf{Step 6}
		
		For every $x\in\R^N$, define
		\[
		u(x)=\sup_{n\in\N}u_n(x).
		\]
		Fix $x\in\R^{N}$.
		If $p(x)=+\infty$, then $u_n(x)=0$ for every $n$, and
		therefore $u(x)=0$. If $p(x)<+\infty$, choose $n_0$ such that $
		|x|\le n_0$ and $ p(x)\le n_0$.
		Then $x\in K_n$ for every $n\ge n_0$, so $
		u_n(x)=0$ for $n\ge n_0$, by \eqref{u_n=0_on_K_n}.
		Consequently,
		\[
		u(x)=\max\{u_1(x),\ldots,u_{n_0-1}(x),0\}<+\infty.
		\]
		Hence the envelope is finite at every point.
		
		\medskip
		\noindent
		\textbf{Step 7}
		
		Fix an arbitrary point $x_0\in\R^{N}$.
		Suppose that $u$ were locally bounded above at $x_0$.
		There would then exist an open neighborhood $W_0$ of $x_0$ and a positive constant $M$ such that
		\begin{equation}\label{u_bounded_above}
				u(x)\le M\qquad(x\in W_0).
		\end{equation}
		Choose a basic ball $B_j$ contained in $W_0$. Because $B_j$ occurs infinitely often in $(V_n)$, choose $n>M$ with $
		V_n=B_j.$
		By construction, $b_n\in B_j\subset W_0$, but 
		\[u(b_n)\ge u_n(b_n)>n>M,\] by \eqref{u_n>n_at_b_n}.
		This contradicts  \eqref{u_bounded_above}.
		Since $x_0$ was arbitrary, $u$ is not locally bounded above at any point of $\R^N$.
	\end{proof}

	\begin{remark}
		Restricting the functions $u_n$ to any nonempty bounded
		open set $\Omega\subset\R^{N}$ gives a pointwise
		finite envelope $u$ that is nowhere locally bounded above
		in $\Omega$. Since the functions $u_n$ are globally defined,
		the counterexample persists even when they are required
		to be defined on a neighborhood of $\overline{\Omega}$.
		
	\end{remark}

	\section{Families of $N-$Separately Subharmonic Functions}
	
	\begin{definition}
		Let $\Omega$ be an open subset of $\C^{N}$. We say that a function
		\[
		v:\Omega\longrightarrow[-\infty,+\infty)
		\]  is \emph{$N$-separately subharmonic}, if the function  $v$ is not identically $-\infty$ on any component of $\Omega$, and for every
			 $j\in\{1,\ldots,N\}$ and every fixed tuple
				\[
				\widehat z_j:=(z_1,\ldots,z_{j-1},z_{j+1},\ldots,z_N)
				\in\C^{N-1},
				\]
				the function
				\[
				\zeta\longmapsto
				v(z_1,\ldots,z_{j-1},\zeta,z_{j+1},\ldots,z_N)
				\]
				is subharmonic, or identically equal to $-\infty$, on every component of the fiber over $\widehat z_j$, namely,
				\[
				\Omega_{\widehat z_j}
				:=
				\bigl\{\zeta\in\C:
				(z_1,\ldots,z_{j-1},\zeta,z_{j+1},\ldots,z_N)
				\in\Omega\bigr\}.
				\]
	\end{definition}
	
	The functions we call $N$-separately subharmonic were considered
	by Lelong in \cite{Lelong}, who called those that are locally
	bounded above functions of class (B). On page 315, he proved
	that every function of class (B) is upper semicontinuous
	(Théorème 1 bis) and subharmonic in the underlying real space $\R^{2n}$ (Théorème 2 bis).

	\begin{theorem}\label{N_sep_sh_Th1}
		Let $\Omega\subset\C^{N}$ be open, and let
		$\{u_\alpha\}_{\alpha\in J}$ be a nonempty family of $N$-separately
		subharmonic functions on $\Omega$. Suppose that
		\[
		u(z):=\sup_{\alpha\in J}u_\alpha(z)<+\infty.
		\] for every $	z\in\Omega$.
		Then the singular set $S$ of $u$ is a relatively closed nowhere dense subset of $\Omega$. Moreover, the upper semicontinuous regularization $u^*$ of $u$ is
		subharmonic on $\Omega\setminus S$, and the set
		\begin{equation}\label{negligible_set_is_polar}
			\{x\in\Omega\setminus S:u(x)<u^*(x)\}
		\end{equation}
		is polar in $\mathbb{R}^{2N}$.
	\end{theorem}
	
	\begin{proof}
		
		To prove that the singular set of $u$ is nowhere dense,
		we may assume that each 
			\begin{equation}\label{u_nonnegative}
			u_\alpha\ge0.
		\end{equation}
		Indeed, replacing
		each $u_\alpha$ by $\max\{u_\alpha,0\}$ preserves
		$N$-separate subharmonicity and replaces the upper envelope
		$u$ by $\max\{u,0\}$, which has exactly the same singular
		set as $u$.

		The proof is by induction on $N$. For $N=1$, the result follows from Corollary \ref*{PlanarUniformBT}.
		Suppose that the theorem is true in $\C^{N-1}$, and we prove it in $\C^{N}$. Let $R$ be the regular set of $u$ (given by Definition \ref{regular_set_singular_set}). This set being open, it  remains only to prove that it is dense. Let $\mathcal{O}$ be an arbitrary nonempty Euclidean open subset
		of $\Omega$. We will prove that 
		\begin{equation}\label{R_is_dense}
			\mathcal{O}\cap R\neq \emptyset.
		\end{equation}
		
		Choose a product
		\[
		D\times W,
		\] relatively compact in $\mathcal{O}$,
		where $D$ is a disc in $\C$ and
		$W$ is a polydisc in $\C^{N-1}$.
		Let $\{K_j\}_{j\in\N}$ be a countable family of compact polydiscs in
		$W$, each having nonempty interior, with the following property:
		every nonempty open subset of $W$ contains at least one $K_j$. For $j,m\geq1$, define
		\begin{align*}
			A_{j,m}:&=
			\bigcap_{w\in K_j}
			\{z_1\in D:u(z_1,w)\leq m\}
		\end{align*}
		We first observe that every $A_{j,m}$ is finely closed in $D$. Indeed, subharmonic functions being finely continuous, $u(\cdot,w)$ is finely lower semicontinuous and each set in the intersection is finely closed, and so is their intersection.
		
		We next prove that
		\begin{equation}\label{A_{j,m}_covers_D}
			D=\bigcup_{j,m\geq1}A_{j,m}.
		\end{equation}
		Fix $z_1\in D$ and  consider the family $
		\{u_\alpha(z_1,\cdot)\}_{\alpha\in J}$ on $W$.
		Every member of this family is $(N-1)$-separately subharmonic and not identically $-\infty$ by \eqref{u_nonnegative}. Moreover,
		for every $w\in W$,
		\[
		\sup_{\alpha\in J}u_\alpha(z_1,w)
		=u(z_1,w)<+\infty,
		\] by our hypothesis.
		The induction hypothesis therefore shows that $
		\{u_\alpha(z_1,\cdot)\}_{\alpha\in J}$ is
		uniformly bounded above on some nonempty open subset $U_{z_1}$ of $W$.
		Thus there exists a constant $C_{z_1}$ such that
		\[
		u_\alpha(z_1,w)\leq C_{z_1}
		\] for every $\alpha\in J$ and every $w\in U_{z_1}$.
		Choose $j$ such that $
		K_j\subset U_{z_1},$
		and choose $m\in\N$ such that $m\geq C_{z_1}$. Then
		$z_1\in A_{j,m}$. This proves (\ref{A_{j,m}_covers_D}).
		
		The fine topology in $\C$ is a Baire topology \cite[Theorem 7.1.3, p. 199]{Armitage Gardiner book}. Since $D$ is not empty and is a
		countable union of the finely closed sets $A_{j,m}$, at least one of these sets, say $A_{j,m}$, has nonempty fine interior. Set \[V=\operatorname{int}_{f}(A_{j,m})\] and choose $z_1^0\in V$.
		Then $V$ is a fine neighborhood of $z_1^0$ and  $V\subset A_{j,m}$.
		By Lemma \ref{Fine-neigh_circl}, there exists a
		Euclidean disc $\Delta$ relatively compact in $ D$, centered at $z_1^0$, such that $
		\partial\Delta\subset V\subset A_{j,m}.$
		It follows from the definition of $A_{j,m}$ that
		\[
		u_\alpha(z_1,w)\leq m
		\]
		for every $	\alpha\in J$ and $(z_1,w)\in\partial\Delta\times K_j.$ Fix $\alpha\in J$ and $w\in K_j$. The
		maximum principle, applied to $u_\alpha(\cdot,w)$ on $\Delta$, gives
		\[
		u_\alpha(z_1,w)\leq m
		\] for all $z_1\in\Delta$.
		Since $\alpha$ and $w$ were arbitrary, we obtain
		\[
		u(z_1,w)\leq m
		\]
		for every
		$(z_1,w)\in\Delta\times K_j$.
		In particular, the family $\{u_\alpha\}_{\alpha\in J}$ is uniformly
		bounded above on the nonempty Euclidean open set $
		\Delta\times\operatorname{int}(K_j).$
		Therefore $
		\Delta\times\operatorname{int}K_j$ is contained in  $R\cap\mathcal{O}$, and \eqref{R_is_dense} follows. The set $R$ is dense in $\Omega$.
		
	To prove the last statement of the theorem, we return to the original
	family $\{u_\alpha\}_{\alpha\in J}$, before the reduction
	\eqref{u_nonnegative}. Recall that the upper semicontinuous regularization
	of $u$ is defined by
	\[
	u^*(z)=
	\limsup_{\substack{\zeta\to z\\(\zeta\in\Omega)}}u(\zeta).
	\]
	Since $u$ is locally bounded above on $\Omega\setminus S$,
	the family $\{u_\alpha\}_{\alpha\in J}$ is locally uniformly
	bounded above there. Consequently, by Lelong's result
	\cite[Théorème~2~bis]{Lelong}, each $u_\alpha$ is subharmonic on
	$\Omega\setminus S$ as a function of $2N$ real variables.
	The fundamental convergence theorem of potential theory
	\cite[Theorem~5.7.1, p.~146]{Armitage Gardiner book}
	then shows that $u^*$ is subharmonic on $\Omega\setminus S$
	and that the negligible set defined in
	\eqref{negligible_set_is_polar} is polar.

	\end{proof}

	\begin{theorem}\label{N_sep_sh_Th2}
		Let $\Omega$ be an open subset of $C^{N}$, and let
		$\{u_\alpha\}_{\alpha\in J}$ be a family of $N$-separately
		subharmonic functions on $\Omega$. Suppose that
		\[
		u(z):=\sup_{\alpha\in J}u_\alpha(z)<+\infty
		\] for every $	z\in\Omega$. If the singular set $S$ of $u$ is not empty, then $S$ is not locally pluripolar at any of its points. In particular, if $S$ is pluripolar, it is empty and $u$ is locally bounded above on $\Omega$.
	\end{theorem}
	
	\begin{proof}
Suppose that $S$ is not empty and is locally pluripolar at some point $\tilde a\in S$. Choose a polydisc
\[
P=D_1\times\cdots\times D_N
\]
relatively compact in $\Omega$ that contains $\tilde a$, and a plurisubharmonic function $\varphi$ on $P$, such that  \(\varphi\not\equiv-\infty\), and
\begin{equation}\label{phi_neg_infty}
	\varphi(z)=-\infty\qquad(z\in S\cap P).
\end{equation}
For every \(a\in S\cap P\) and \(1\le j\le N\), define the $j-$th slice	
	\[
	S_j(a)=\{\zeta\in D_j:(a_1,\ldots,a_{j-1},\zeta,a_{j+1},\ldots,a_N)\in S\}.
	\]

\textbf{Claim}: Each $S_j(a)$ is nonpolar, for every $a\in S\cap P$.

Suppose $S_j(a)$ is polar for some $a\in S$. A polar set being thin at every of its points, we can use Lemma \ref{Fine-neigh_circl} to choose a disc $\Delta$ relatively compact in $ D_j$ and centered at \(a_j\) such that
\[\partial\Delta\cap S_j(a)=\varnothing.\]
Let $a':=(a_1,\dots,a_{j-1})	$ and $a'':=(a_{j+1},\dots,a_N)$. That means the compact set $\{a'\}\times\partial\Delta\times\{a''\}$ lies in the open set $\Omega\setminus S$, and so there are neighborhoods $W_1$ and $W_2$ of $a'$ and $a''$, respectively, such that \[K:=\overline{W_1}\times\partial\Delta\times\overline{ W_2}\] is contained in $\Omega\setminus S$. Since $K$ is compact, there exists a constant $M$ satisfying 
$$u(z)\le M\qquad(z\in K).$$
	Fix $z'=(z_1,\dots,z_{j-1})$ in $ W_1$ and $z''=(z_{j+1},\dots,z_n)$ in $W_2$. Fix also $\alpha\in J$. Applying the maximum principle to the subharmonic function $u_\alpha(z',\cdot,z'')$, then taking the supremum over $\alpha$,  shows that the above inequality holds on $W_1\times \Delta\times W_2$, which is a neighborhood of $a$. This contradicts the choice of $a\in S$, and proves the claim.

To prove the theorem, we will propagate the identity \(\varphi=-\infty\) coordinate by coordinate using induction. Choose \(b=(b_1,\ldots,b_N)\in P\) with $
	\varphi(b)>-\infty.$ We shall prove by induction over $k$ that, for every \(a=(a_1,\ldots,a_N)\in S\cap P\), and every $k=0,\ldots,N$, 
	\[
	\varphi(b_1,\ldots,b_k,a_{k+1},\ldots,a_N)=-\infty.
	\]
	
The case of \(k=0\) follows from \eqref{phi_neg_infty}.
Fix \(0\le k<N\), and suppose the assertion holds at step \(k\) for every \(a\in S\cap P\). Fix $a\in S\cap P$ and consider the function
	\[
	v(\zeta)=
	\varphi(b_1,\ldots,b_k,\zeta,a_{k+2},\ldots,a_N).
	\] For any $\zeta\in S_{k+1}(a)$, the point 
	\[
	a'=(a_1,\ldots,a_k,\zeta,a_{k+2},\ldots,a_N)
	\] belongs to $S\cap P$. Applying the induction hypothesis to $a'$ gives
	\[
	\varphi(b_1,\ldots,b_k,\zeta,a_{k+2},\ldots,a_N)
	=-\infty.
	\]
	But this is precisely
	\[
	v(\zeta)=-\infty
	\] on $S_{k+1}(a)$. Since $\varphi$ is plurisubharmonic, $v$ is subharmonic or identically $-\infty$. Since $S_{k+1}$ is not polar, according to the Claim, it follows that $v\equiv-\infty$, and in particular,
	\[
	\varphi(b_1,\ldots,b_{k+1},a_{k+2},\ldots,a_N)=-\infty.
	\] This completes the induction.
	
	Taking \(k=N\), we obtain $
	\varphi(b_1,\ldots,b_N)=-\infty,$ contradicting the choice of $b$. Therefore $S\not=\varnothing$ implies that $S$ is not locally pluripolar at $\tilde{a}$, and since $\tilde{a}$ was arbitrary, at any point of $S$.  
	\end{proof}

		\begin{corollary}
		Let $\Omega$ be an open subset of $\C^{N}$, and let
		$\{u_\alpha\}_{\alpha\in J}$ be a family of 	plurisubharmonic functions on $\Omega$. Suppose that
		\[
		u(z):=\sup_{\alpha\in J}u_\alpha(z)<+\infty
		\] for every $	z\in\Omega$. If the singular set $S$ of $u$ is not empty, then $S$ is not locally pluripolar at any of its points. In particular, if $S$ is pluripolar, it is empty and $u$ is locally bounded above on $\Omega$.
	\end{corollary}
	\begin{proof}
		Every plurisubharmonic function is \(N\)-separately subharmonic, so the conclusion follows from Theorem \ref{N_sep_sh_Th2}.
	\end{proof}

	\begin{corollary}
		Let $\Omega\subset\C^{N}$ be open, and let
	$u$ be an $N$-separately
	subharmonic function on $\Omega$. Then the singular set $S$ of $u$ is relatively closed and nowhere dense in $\Omega$. Moreover,  if $S$ is pluripolar, it is empty and $u$ is subharmonic on $\Omega$.
	\end{corollary}
	\begin{proof}
		Apply Theorems \ref{N_sep_sh_Th1} and \ref{N_sep_sh_Th2} to the singleton family $\{u\}$. If $S$ is empty, then $u$ is locally bounded above on $\Omega$, and hence subharmonic by Lelong's theorem.
	\end{proof}
	\begin{remark}
			Let \(P=D_1\times\cdots\times D_N\) be a polydisc in $\C^N$, and let \(u\) be \(N\)-separately subharmonic on an open neighborhood of the closure \(\overline P\) of $P$. Let  $\Gamma=\partial D_1\times\cdots\times\partial D_N$ be the distinguished boundary of $P$.
		Then
		\[\sup_{\zeta\in \overline P}u(\zeta)=\sup_{\zeta\in\Gamma}u(\zeta)\]
	where the common supremum may be \(+\infty\). In particular, if \(u\) is bounded above on \(\Gamma\), then \(u\) is subharmonic on \(P\).
	\end{remark}

	\begin{proof}
	 Applying the one-variable maximum principle successively in each coordinate gives, for every \(z\in\overline P\),
	\[
	\begin{aligned}
		u(z)
		&\le \sup_{\zeta_1\in\partial D_1}
		u(\zeta_1,z_2,\ldots,z_N)\\
		&\le \sup_{\substack{\zeta_1\in\partial D_1\\
				\zeta_2\in\partial D_2}}
		u(\zeta_1,\zeta_2,z_3,\ldots,z_N)\\
		&\le\cdots\le\sup_\Gamma u.
	\end{aligned}
	\]
	The reverse inequality follows from \(\Gamma\subset\overline P\). If the boundary supremum is finite, this bounds \(u\) above throughout \(P\), so Lelong’s theorem gives joint subharmonicity there.
	\end{proof}

\begin{example}
	Using Wiegerinck--Zeinstra's construction \cite{Wiegerinck Zeinstra},
	we give, for every $N\ge2$, an example of an $N$-separately
	subharmonic function that is not jointly subharmonic.
	Take their nonnegative subharmonic functions $v_n$ on $\C$,
	which vanish eventually at each fixed point and satisfy
	$v_n(a_n)\to+\infty$ for some sequence $a_n\to0$. Then
	\[
	u(z_1,\ldots,z_N):=\sum_{n=1}^{\infty}\prod_{j=1}^{N}v_n(z_j)
	\]
	is finite everywhere and subharmonic in each variable.
	Indeed, fixing all coordinates except $z_j$, we obtain a finite
	sum of nonnegative multiples of the functions $v_n(z_j)$,
	because $N\ge2$ and $v_n(z_k)=0$ for all sufficiently large $n$
	at any fixed coordinate $z_k$ with $k\ne j$.
	However,
	\[
	u(a_n,\ldots,a_n)\ge \bigl(v_n(a_n)\bigr)^N\longrightarrow+\infty,
	\]
	so $u$ is not locally bounded above at the origin and hence
	is not jointly subharmonic on any neighborhood of the origin.
\end{example}
	
	\section{Characterization of Singular Sets}

	\subsection{Families Without a Continuity Assumption}
	
	\begin{definition}\label{admissible_def}
	Let \(X\) be a locally compact, second-countable Hausdorff space. A class \(\mathcal E\) of functions from $X$ to $[-\infty,+\infty)$ is called \textit{admissible} if the following holds:
	\begin{itemize}
	\item [(i)]	\(0\in\mathcal E\);
	\item [(ii)]	every member of \(\mathcal E\) is locally bounded above;
	\item [(iii)]	\(\varphi\in\mathcal E\Rightarrow a\varphi+c\in\mathcal E\) for \(a>0\), \(c\in\mathbb R\);
	\item [(iv)]	\(\varphi\in\mathcal E\Rightarrow\max\{\varphi,0\}\in\mathcal E\).
		
	\end{itemize}
	\end{definition}

Let $\Omega$ be an open subset of $\mathbb R^N$ or $\mathbb C^N$.
The class of upper semicontinuous functions defined on  $\Omega$ with values in $[-\infty,+\infty)$
is admissible. Such functions are locally bounded above, and the
operations in Definition~\ref*{admissible_def} preserve upper
semicontinuity. In the real case, $\mathcal{SH}(\Omega)$ is also
admissible, since these operations preserve subharmonicity. In the
complex case, the same reasoning applies to $\mathcal{PSH}(\Omega)$.
	
	Let $X$ and $\mathcal{E}$ be as in the above definition. Recall that a subset $A$ of $X$ is relatively compact if its closure $\overline{A}$ is compact in $X$. For a nonempty relatively compact set \(A\subset X\), we define the \textit{extended hull }of $A$ with respect to $\mathcal{E}$ by
	\begin{equation}\label{def_hull_for_bounded A}
			\widehat A_{\mathcal E}
		=	\left\{x\in X:\varphi(x)\le\sup_A \varphi
		\text{ for every }\varphi\in\mathcal E\right\},
		\end{equation}
			with \(\widehat\varnothing_{\mathcal E}=\varnothing\).
	\begin{definition}\label{def_extended_hull_sep}
		A closed set $S$ in $X$ has the \textit{extended hull separation} property with respect to \(\mathcal E\) if there exist increasing sequences of relatively compact sets $(T_n)_{n\in\N}$ in $ S $ and compact sets $(C_n)_{n\in\N}$ in $ X\setminus S$ such that:
		\begin{itemize}
			\item [(i)]$S=\bigcup_{n\in\N}T_n$;
			\item [(ii)]every compact subset of \(X\setminus S\) is contained in \(C_n\) for all sufficiently large \(n\);
			\item [(iii)]$S\cap\operatorname{int}_X\left(\left(\widehat{T_n\cup C_n}\right)_{\mathcal E}\right)=\varnothing$ for every $n\in\N$.
		\end{itemize}
		
	\end{definition}
	
	\begin{theorem}\label{Theorem_for_mathcal_E}
			Let $X$ be a locally compact, second-countable Hausdorff space.  Suppose $\mathcal{E}$ is an admissible class of functions on $X$, and let $S$ be a closed subset of $X$. The following are equivalent:
		
		(i) There exists a nonempty family $\{u_\alpha\}_{\alpha\in J}$ of functions of $\mathcal E$ such that
		\begin{equation}\label{u_less_than_infinity}
				u(x):=\sup_{\alpha\in J}u_\alpha(x)<+\infty
		\end{equation}
		for all $x\in X$ and $S$ is the singular set of $u$.
		
		(ii) The set $S$ has the extended hull separation property with respect to $\mathcal E$.
		
		Moreover, the family in (i) may be chosen to be a sequence of nonnegative functions.
	\end{theorem}
	\begin{proof}
		Let $(L_n)_{n\in\N}$ be an increasing compact exhaustion of $X$ whose interiors cover $X$.	If \(S=\varnothing\), condition (i) of the Theorem holds with \(u_n\equiv0\), and condition (ii) holds with $
		T_n=\emptyset$ and $C_n=L_n$.
		We therefore assume \(S\ne\varnothing\).
		
		\((i)\Rightarrow(ii)\). Since \(X\setminus S\) is an open subset of a locally compact, second-countable Hausdorff space, choose an increasing sequence of compact sets \(C_n\subset X\setminus S\) whose interiors relative to \(X\setminus S\) cover \(X\setminus S\). Every compact subset of \(X\setminus S\) is then contained in \(C_n\) for all sufficiently large \(n\).	If \(X\setminus S=\varnothing\), take \(C_n=\varnothing\). This proves condition (ii) of Definition ~ \ref{def_extended_hull_sep}.
		
		By hypothesis, \(u\) is locally bounded above at every point of \(X\setminus S\). Compactness therefore implies that \(u\) is bounded above on each nonempty compact \(C_n\). Therefore, by defining  
		\[M_n:=\max  \left\{ \sup_{x\in C_n}u(x),n\right\},\] with the convention \(\sup\varnothing=-\infty\), 
		each \(M_n\) is finite. Define also,
			\begin{equation}\nonumber
			T_n=S\cap L_n\cap\{x\in X:  u(x)\le M_n\}.
		\end{equation}
		Since the sequence $(C_n)_{n\in\N}$ is increasing, so is the sequence $(M_n)_{n\in \N}$. Since the sequences $(M_n)_{n\in \N}$ and $(L_n)_{n\in\N}$ are increasing, so is the sequence $(T_n)_{n\in\N}$.	Since $T_n$ is contained in the compact set $ S\cap L_n$, the closure $\overline{T_n}$ is a compact subset of $S$ and $T_n$ is relatively compact in $S$.
			
			 Let  \(x\in S\). We have \(x\in L_n\) for all sufficiently large \(n\). Since \(u(x)<+\infty\) by \eqref{u_less_than_infinity} and \(M_n\ge n\), we also have \(u(x)\le M_n\) eventually. Hence \(x\in T_n\) for all sufficiently large \(n\) and 
		\[
		\bigcup_{n\ge1}T_n=S,
		\]
		which proves condition (i) of Definition \ref*{def_extended_hull_sep}.

		Fix $n\in\N$, and write \(A_n=T_n\cup C_n\). If \(A_n=\varnothing\), condition (iii) of Definition \ref{def_extended_hull_sep} follows from our convention for the hull of empty set. Otherwise, for every \(\alpha\in J\), $u_\alpha\le u\le M_n$ on $A_n$, by definitions of $M_n$ and $T_n$.
		Consequently, it follows from the definition of the hull that
		\[	u_\alpha(x)\le\sup_{A_n}u_\alpha\le M_n \]
        for $x\in \widehat{A_n}_{\mathcal E}.$
		Taking the supremum over \(\alpha\) gives
		\[
		u(x)\le M_n \qquad(x\in\widehat{A_n}_{\mathcal E})
		.\]
	This inequality proves that the interior of $\widehat{A_n}_{\mathcal E}$ does not intersect $S$: if a point of \(S\) belonged to the interior in $X$ of this hull, it would have a neighborhood on which \(u\) is bounded above, contradicting the definition of \(S\). Thus
		\[
		S\cap\operatorname{int}_X
	\left(\left(	\widehat{T_n\cup C_n}\right)_{\mathcal E}\right)
		=\varnothing.
		\]
		This proves condition (iii) of Definition \ref{def_extended_hull_sep}.
		
		\vspace*{5mm}
		\((ii)\Rightarrow(i)\). Let $(T_n)_{n\in\N}$ and $(C_n)_{n\in\N}$ satisfy the extended hull separation property, and set
		\[
		A_n=T_n\cup C_n.
		\]
		Choose a countable basis of nonempty open sets for $X$. Enumerate those basic open sets that meet $S$ as a sequence $(V_n)_{n\in\N}$ in which each such basic open set occurs infinitely often.
		
		For every \(n\),
		\[
		V_n\setminus\widehat{A_n}_{\mathcal E}\ne\varnothing.
		\]
		Indeed, if \(V_n\subset\widehat {A_n}_\mathcal{E}\), then, because \(V_n\) is open, $
		V_n\subset\operatorname{int}_X\widehat{A_n}_{\mathcal E}$.
		As \(V_n\cap S\ne\varnothing\), this contradicts condition (iii). We may therefore choose
		\[
		b_n\in V_n\setminus\widehat{ A_n}_{\mathcal E}.
		\]
		
		Now, following Wiegerinck's construction \cite{Wiegerinck}, we construct a nonnegative function $u_n$ in $\mathcal E$ that vanishes on $A_n$ and satisfies $u_n(b_n)=n+1$. Suppose first that \(A_n\ne\varnothing\). By the definition of the hull, there exists \(v_n\in\mathcal E\) such that $v_n(b_n)>\sup_{A_n}v_n$.
		In particular, \(v_n(b_n)\) is a finite number. Choose a real number $c_n$ satisfying
		\begin{equation}\nonumber
		\sup_{A_n}v_n<c_n<v_n(b_n).
		\end{equation}
		Such a choice is possible even if \(\sup_{A_n}v_n=-\infty\). Define
		\[
		u_n(x)=
		\frac{n+1}{v_n(b_n)-c_n}
		\max\{v_n(x)-c_n,0\}.
		\]
	By admissibility of $\mathcal E$, we have  $u_n\in\mathcal E$. Furthermore, $u_n$ is nonnegative and finite-valued, and
		\begin{equation}\label{u_n(b_n_=n+1)}
		u_n=0\quad\text{on }A_n,
		\qquad
		u_n(b_n)=n+1,
		\end{equation}as desired.
		If \(A_n=\varnothing\), define \(u_n\equiv n+1\); these same assertions hold, with the vanishing assertion being vacuous.
		
		Set
		\[
		u(x):=\sup_{n\in\N}u_n(x).
		\]
		We verify that \(u\) is finite everywhere and has singular set \(S\).

		Let \(x\in X\setminus S\). Since $S$ is closed, by local compactness and the Hausdorff property, choose a compact neighborhood $K$ of $x$ contained in $X\setminus S$.
		By condition (ii), there exists \(n_0\) such that $ C_n$ contains the compact $K$ for every $n\ge n_0$. Thus \(u_n=0\) on \(K\) for \(n\ge n_0\), and
		\[
		u=\max\{u_1,\ldots,u_{n_0-1},0\}
		\] on $K$.
		Each of the finitely many functions on the right is bounded above on $K$, since it belongs to $\mathcal{E}$ and by compactness of $K$. Therefore \(u\) is bounded above on $K$ and so locally bounded above at \(x\).
		
		Next, fix \(x\in X\). If \(x\in S\), condition (i) and the increasing nature of $(T_n)_{n\in\N}$ imply that \(x\in T_n\subset A_n\) for all sufficiently large \(n\). Thus $u_n(x)=0$ for all sufficiently large $n$. 
		Since each \(u_n(x)\) is finite, it follows that \(u(x)<+\infty\). If $x\in X\setminus S$,  the previous paragraph already shows that $u(x)<+\infty$.
		
		Finally, let \(x\in S\), and let \(U\) be any open neighborhood of \(x\). Choose a basic open set \(V\) with $x\in V\subset U$.
		This open set meets $S$, so it occurs infinitely often in \((V_n)\), at arbitrarily large indices. For every such index, $b_n\in V_n=V$, which is contained in $U$, and 
		\[u(b_n)\ge u_n(b_n)=n+1,\] by \eqref{u_n(b_n_=n+1)}.
		Hence \(u\) is unbounded above on \(U\). Since \(U\) was arbitrary, \(x\) belongs to the singular set of \(u\). Thus the singular set of \(u\) is exactly \(S\). Moreover, the construction provides the asserted sequence of nonnegative functions.
	\end{proof}

	\begin{corollary}\label{coroll_for_sh}
		Let $\Omega$ be an open subset of $\R^N$ with $N\ge2$, or of $\C^N$ with $N\ge1$. Set $\mathcal E=\mathcal{SH}(\Omega)$ in the real case and $\mathcal E=\mathcal{PSH}(\Omega)$ in the complex case. For a relatively closed set $S\subset\Omega$, the following are equivalent.
		\begin{itemize}
			\item [(i)]There exists a sequence $(u_n)_{n\in\N}$ of nonnegative functions of $\mathcal E$ such that
			\begin{equation}\nonumber
				u(\zeta):=\sup_{n\ge1}u_n(\zeta)<+\infty
			\end{equation}
			for all $\zeta\in\Omega$, and $S$ is the singular set of $u$.
			\item [(ii)] The set $S$ has
			the extended hull separation property with respect to
			$\mathcal E$.
		\end{itemize}
		
	\end{corollary}
	\begin{proof}
		The space $\Omega$ is locally compact, second-countable, and Hausdorff, and $\mathcal E$ is admissible. Apply Theorem \ref{Theorem_for_mathcal_E}.
	\end{proof}

	\subsection{Families of Continuous Functions}
We now characterize the singular sets of upper envelopes for admissible classes of continuous functions.
	
	\begin{definition}\label{hull_separation_def}
		Let $X$ be a locally compact, second-countable Hausdorff space, and let  $\mathcal F$ be an admissible class of continuous real-valued functions on $X$. A closed set $S$ in $ X$ has the \textit{hull separation property} with respect to \(\mathcal F\) if there exist increasing sequences \((T_n)_{n\in\N}\) and \((C_n)_{n\in\N}\) of compact subsets of \(S\) and \(X\setminus S\), respectively, such that: 
		\begin{itemize}
			\item[(i)] $\bigcup^\infty_{n=1}T_n=S$; 
			
			\item [(ii)] every compact
			subset of $X\setminus S$ is contained in $C_n$
			for all sufficiently large $n$;
			\item[(iii)] for every $n\in\N$
		$S\cap\operatorname{int}_X\left(\left(\widehat{T_n\cup C_n}\right)_{\mathcal F}\right)=\varnothing$ for every $n\in\N$.
		\end{itemize}
	\end{definition}
	Here the hull is defined by \eqref{def_hull_for_bounded A} with \(\mathcal E=\mathcal F\) and \(\widehat\varnothing_{\mathcal F}=\varnothing\).
	
	It is well known that if $X$ is a Baire space, each $u_\alpha$ is continuous and $u=\sup_{\alpha\in J}u_\alpha$ is finite everywhere, then the singular set of $u$ is closed and nowhere dense in $\Omega$. This follows from the Baire category theorem and in fact, lower semicontinuity of the functions $u_\alpha$ suffices (see Bourbaki \cite[ Theorem 2, p. 194]{Bourbaki}). The following result provides, in a certain sense, a converse to this result.
	
	\begin{theorem}\label{S implies HSP}
		Let $X$ be a locally compact, second-countable Hausdorff space, and let $\mathcal F$ be an admissible class of continuous real-valued functions on $X$.  Let $S$ be a closed subset of $X$. The following are equivalent:
			\begin{itemize}
			\item [(i)] There exists a nonempty family $\{u_\alpha\}_{\alpha\in J}$ of functions of $\mathcal F$ such that
		\[
		u(x):=\sup_{\alpha\in J}u_\alpha(x)<+\infty
		\] for all $x\in X$, and $S$ is the singular set of $u$.
		
		\item [(ii)] The set $S$ has
	the hull separation property with respect to $\mathcal{F}$.
		\end{itemize}
		Moreover, the family in (i) may be chosen to be a sequence of nonnegative functions.
		
	\end{theorem}
	
	\begin{proof}
		The proof of Theorem \ref*{Theorem_for_mathcal_E}  applies with minor modifications.
	For \((i)\Rightarrow(ii)\), use the necessity proof of Theorem 6.3 with \(\mathcal E=\mathcal F\). Its sets
	\[
	T_n=S\cap L_n\cap\{u\le M_n\}
	\]
	are compact: indeed, the last set on the right is closed, since by continuity of each $u_\alpha$, the upper envelope $u$ is lower semicontinuous. Thus the construction satisfies Definition \ref{hull_separation_def}.
	
	For $(ii)\Rightarrow(i)$, the hull separation property implies
	the extended hull separation property, since compact sets are
	relatively compact. Theorem~6.3, applied to $\mathcal F$, therefore
	gives the required sequence of nonnegative functions $u_n$.
	Each $u_n$ is continuous, as it belongs to $\mathcal F$.
		
	\end{proof}

	\begin{remark}
		
		In Corollary \ref{coroll_for_sh}, let $\mathcal F$ be the class of functions in $\mathcal E$ that are continuous on $\Omega$. By Theorem \ref{S implies HSP}, a relatively closed set $S$ of $\Omega$ is the singular set of a pointwise bounded-above envelope of a sequence of nonnegative functions in $\mathcal F$ if and only if $S$ has the hull separation property with respect to $\mathcal F$.
		
	\end{remark}

\subsection{Sufficient Conditions}
	
	So far we have characterized the singular set for families of classes $\mathcal{E}$ and $\mathcal{F}$. The following result gives a sufficient condition for a set to be the singular set of the upper envelope of a family of subharmonic functions. 
	
	\begin{theorem}\label{sufficient cond for S}
		Let $N\ge3$, and let $S$ be a closed subset of
		$\R^{N}$. Suppose that there exists a polar set
		$D\subset S$ which is dense in $S$ and whose
		components are all unbounded. Then there exists a sequence $(u_n)_{n\in\N}$ of nonnegative subharmonic functions on $\R^{N}$
		such that
		\[
		u(x):=\sup_{n\ge1}u_n(x)<+\infty
		\]
		for all $x\in\R^N$, and the singular set of $u$ is exactly $S$. 
	\end{theorem}
	
	\begin{proof}
		If $S=\varnothing$, take $u_n\equiv0$. We therefore assume
		that $S\ne\varnothing$. Since $D$ is polar, choose a nonnegative
		superharmonic function $p$ on $\R^N$ such that $p\not\equiv+\infty$ and $ p=+\infty$ on $D$. For $n\ge1$, define
		\[
		A_n=\{x\in\R^N:  |x|\le n,\ p(x)\le n\},
		\]
		\[
		C_n=
		\left\{x\in\R^N: |x|\le n,\
		\operatorname{dist}(x,S)\ge\frac1n\right\},
		\] and 
		\[	K_n=A_n\cup C_n.\]
		These sets are compact and increase with $n$.
		Compactness of $A_n$ follows from lower semicontinuity
		of $p$.
		
		Notice that every point where $p$ is finite eventually belongs to $A_n$ and every compact subset of
		$\R^{N}\setminus S$ eventually belongs to $C_n$: the reason is that 
		such a compact set is bounded and has positive distance from the closed set $S$. Also,
		since $D\subset S$ and $p=+\infty$ on $D$,
		\begin{equation}\label{D_cap_K_n_empty}
			D\cap K_n=\varnothing.
		\end{equation}
		Let $\widehat{K}_n$ be the filled hull of $K_n$; that is the union of $K_n$ and all bounded components of $\R^{N}\setminus K_n$. The set $\widehat{K}_n$ is compact and its complement is connected.
		
		\medskip
		
		We claim that
		\begin{equation}\label{D_cap_F_n_nonempty}
			D\cap \widehat{K}_n=\varnothing.
		\end{equation}
		Indeed, let $E$ be any component of $D$.
		Since $E$ is connected and disjoint from $K_n$ by (\ref{D_cap_K_n_empty}), it lies
		in a single component of
		$\R^{N}\setminus K_n$. That component must be unbounded because $E$ is unbounded. Thus $E$ is not contained in any of the bounded components adjoined
		to form $\widehat{K}_n$. This and (\ref{D_cap_K_n_empty}) prove the claim. Also, since $D$ is dense in $S$,
		\begin{equation}\label{S_int_hat_K_n_empty}
			S\cap\operatorname{int}(\widehat{K}_n)=\varnothing.
		\end{equation}
		For if $a\in S\cap\operatorname{int}(\widehat{K}_n)$, an open ball
		about $a$ contained in $\widehat{K}_n$ would meet $D$, contradicting
		$D\cap \widehat{K}_n=\varnothing$.

		Now choose a countable basis of open balls in $\R^{N}$. Enumerate the basic balls meeting $S$ as a sequence $(V_n)_{n\in\N}$ in which each such ball occurs infinitely often.
		For every $n\in\mathbb{N}$, the open set $V_n\setminus\widehat{K}_n$ is nonempty. Indeed, otherwise we would have $V_n\subset\operatorname{int}(\widehat{K}_n)$, since $V_n$ is open. As $V_n$ meets $S$, this would contradict (\ref{S_int_hat_K_n_empty}). Since $p$ is not identically $+\infty$ on any nonempty open set, we may therefore choose $
		b_n\in V_n\setminus\widehat{K}_n
		$ with $p(b_n)<+\infty$.
		
		Moreover, $\R^{N}\setminus\widehat{K}_n$ is connected and open. Since $N\ge2$, removing a single point preserves its connectedness. So
		\[
		\R^{N}\setminus\bigl(\widehat{K}_n\cup\{b_n\}\bigr)
		\]
		is also connected and open.
		
		As in Step 4 of Theorem \ref{S_Can_be_R}, one proves that there exists a harmonic function $h_n$ satisfying $|h_n|\le1/2$ on $\widehat{K}_n$ and $h_n(b_n)>p(b_n)+n+3/2$. Define $u_n$ as in Step~5 of Theorem~4.1, with $K_n$
		replaced by $\widehat K_n$. The argument in Step~6
		shows that $u=\sup_n u_n$ is finite everywhere.
		Every compact subset of $\mathbb R^N\setminus S$
		is eventually contained in $C_n\subset\widehat K_n$,
		where $u_n=0$. Thus $u$ is locally bounded above
		outside $S$. The argument in Step~7, applied to
		the basic balls meeting $S$, shows that $u$ is
		unbounded above on every neighborhood of each
		point of $S$. Hence its singular set is exactly $S$.
		
	\end{proof}


	\section{Separately Subharmonic Functions}

	We apply the results on the families of subharmonic functions to the class of separately subharmonic functions. 
	\begin{definition}
		Let $N_1,N_2\ge2$ be natural numbers and $\Omega_1$ and  $\Omega_2$ be open sets in $\R^{N_1}$ and $\R^{N_2}$, respectively. A function 
		$$u:\Omega_1\times\Omega_2\to[-\infty,+\infty)$$
		is called \textit{separately subharmonic} if for each fixed $(a,b)\in \Omega_1\times\Omega_2$, the partial functions $u(\cdot,b)$ and $u(a,\cdot)$ are subharmonic on $\Omega_1$ and $\Omega_2$, respectively.
	\end{definition}
	
	We will use the following result due to Avanissian \cite{Avanissian}: \textit{A separately subharmonic function defined on an open set $\Omega_1\times\Omega_2$ of  $\R^{N_1}\times\R^{N_2}$, with $N_1,N_2\ge2$, that is locally bounded above on  $\Omega_1\times\Omega_2$ is subharmonic.}
	
	\begin{theorem}\label{Th 7_2}
		Let $\Omega_1$ and
		$\Omega_2$ be nonempty bounded open sets in $\C$ and $\R^N$, respectively, with $N\geq2$.
		Let $u$ be separately subharmonic on a neighborhood of
		$\overline{\Omega_1\times\Omega_2}$.
		Then, there exists a relatively closed nowhere dense set $S_1$ in $\Omega_1$ such that $u$ is subharmonic on $(\Omega_1\setminus S_1)\times\Omega_2$. Moreover, if $S_1$ is polar, then $u$ is subharmonic on $\Omega_1\times\Omega_2$.
		
	\end{theorem}
	
	\begin{proof}
	 For $z\in\Omega_1$, define
			\[
		\Phi(z):=
		\max_{y\in\overline{\Omega_2}}u(z,y).
		\]
		For fixed $z\in\Omega_1$, by upper semicontinuity of $u(z,\cdot)$, we have 
		\begin{equation}\label{Phi_finite}
			\Phi(z)<+\infty.
		\end{equation}
		By Corollary \ref{PlanarUniformBT}, applied to the family $\{u(\cdot,y)\}_{y\in\overline{\Omega_2}}$ of subharmonic functions, the singular set \(S_1\) of \(\Phi\) is relatively closed and nowhere dense in \(\Omega_1\). Since \(u(z,y)\le\Phi(z)\) for all \((z,y)\in\Omega_1\times\Omega_2\), the function \(u\) is locally bounded above on \((\Omega_1\setminus S_1)\times\Omega_2\). Avanissian’s theorem therefore gives subharmonicity there.

\smallskip

\smallskip

		  Suppose $S_1$ is polar. Fix $(a,b)\in
		\Omega_1\times\Omega_2$. We will prove that $u$ is locally bounded above at $(a,b)$.	Let $V:=(\Omega_1\setminus S_1)\cup\{a\}$. Notice that $V=\Omega_1\setminus(S_1\setminus\{a\})$ is a fine neighborhood of $a$. Indeed, $S_1\setminus\{a\}$ being polar, is finely closed. So $V$ is finely open and contains $a$.  	It follows from Lemma \ref{Fine-neigh_circl} that there exists a Euclidean disc $\Delta_1$, centered at $a$ and relatively compact in $\Omega_1$, whose boundary $\partial\Delta_1$ lies in $V$. Since $a$ is not in $\partial\Delta_1,$ we obtain
		\begin{equation}\label{boundary_delta_regular}\nonumber
			\partial\Delta_1\cap S_1=\emptyset.
		\end{equation}
		
		Thus $\Phi$ is locally bounded above at each point of the compact $\partial \Delta_1$. By compactness there exists a constant $M$ such that 
		\[u(z,y)\le \Phi(z)\le M\] for all $z\in \partial\Delta_1$ and $y\in\Omega_2$. Fix $y\in\Omega_2$. The above inequality and the maximum principle applied to the subharmonic function $u(\cdot,y)$ imply that
	$u(z,y)\le M$ for all $z\in \Delta_1$. Since $y$ was arbitrary, the last inequality holds in $\Delta_1\times\Omega_2$, which is a neighborhood of $(a,b)$. Thus $u$ is locally bounded above at any arbitrary point of $\Omega_1\times\Omega_2$, and  $u$ is subharmonic everywhere, by Avanissian's theorem.

	\end{proof}

	\begin{corollary}
		Let $\Omega_1,\Omega_2$ be bounded open sets in $\C$,
		and let $u$ be separately subharmonic on an open neighborhood
		of $\overline{\Omega_1\times\Omega_2}$.
		There exist relatively closed nowhere dense sets
		$S_1$ in $\Omega_1$ and $S_2$ in $\Omega_2$ such that $u$
		is subharmonic on
		\[
		(\Omega_1\times\Omega_2)\setminus(S_1\times S_2).
		\]
	\end{corollary}

	\begin{proof}
	Apply Theorem \ref{Th 7_2} to $u$ and then apply it again with the variables interchanged. This gives subharmonicity on each of the two open sets on the left below. Since subharmonicity is local, $u$ is subharmonic on their union:
	\[
	\bigl((\Omega_1\setminus S_1)\times\Omega_2\bigr)
	\cup
	\bigl(\Omega_1\times(\Omega_2\setminus S_2)\bigr)
	=
	(\Omega_1\times\Omega_2)\setminus(S_1\times S_2).
	\]
		
	\end{proof}

	The following counterexample  proves that the previous theorem and its corollary do not hold for $N_1,N_2\ge3$. We adapt the construction from the proof of Theorem \ref{S_Can_be_R}, using a series of products as in the counterexample of Wiegerinck \cite{Wiegerinck} (see also Wiegerinck and Zeinstra \cite{Wiegerinck Zeinstra}).

	\begin{theorem}\label{Th 7_4}
		For every pair of integers $N_1,N_2\ge3$, there exists a
		nonnegative separately subharmonic function $u$ on
		$\R^{N_1}\times\R^{N_2}$ which is not
		subharmonic at any point of
		$\R^{N_1}\times\R^{N_2}$.
	\end{theorem}

	\begin{proof}
		We use the construction in the proof of
		Theorem~\ref{S_Can_be_R}, with superscript $(\nu)$
		for the corresponding objects in $\R^{N_\nu}$,
		$\nu=1,2$.
		
		Choose the sequences $(V_n^{(1)})_{n\in\N}$ and
		$(V_n^{(2)})_{n\in\N}$ simultaneously so that
		$(V_n^{(1)}\times V_n^{(2)})_n$ enumerates a
		countable basis of open rectangles in
		$\R^{N_1}\times\R^{N_2}$,
		with each basic rectangle occurring infinitely often.
		Using these sequences in the construction, we obtain by Step 5
		nonnegative subharmonic functions $u_n^{(\nu)}$ and
		points $b_n^{(\nu)}\in V_n^{(\nu)}$ such that
		\begin{equation}\label{u_n(b_n) biger than n}
			u_n^{(\nu)}(b_n^{(\nu)})>n.
		\end{equation}
		Moreover, by the same Step 5, for each fixed $x^{(\nu)}\in\R^{N_\nu}$,
		we have 
		\begin{equation}\label{u_nu_(x_nu)_zero_n_big}
			u_n^{(\nu)}(x^{(\nu)})=0
		\end{equation} for all sufficiently large $n$.
		
		Define
		\[
		u(x^{(1)},x^{(2)})
		=\sum_{n=1}^{\infty}
		u_n^{(1)}(x^{(1)})u_n^{(2)}(x^{(2)}).
		\]
		Whenever either variable is fixed, this series reduces
		to a finite sum, by (\ref{u_nu_(x_nu)_zero_n_big}). Thus $u$ is finite-valued and
		separately subharmonic, since each factor in the summand is a finite
		linear combination of subharmonic functions with
		nonnegative coefficients.
		
		Suppose that $u$ is locally bounded above at some
		point $(a^{(1)},a^{(2)})$. Then there exist open
		neighborhoods $\omega^{(\nu)}$ of $a^{(\nu)}$ and
		$M>0$ such that $
		u\le M$ on $\omega^{(1)}\times\omega^{(2)}$.
		By our choice of the sequences $(V_n^{(\nu)})_{n\in\N}$,
		there exists $n>M$ such that
		\[
		(b_n^{(1)},b_n^{(2)})\in V_n^{(1)}\times V_n^{(2)}\subset\omega^{(1)}\times \omega^{(2)}.\]
		Hence, by \eqref{u_n(b_n) biger than n},
		\[
		u(b_n^{(1)},b_n^{(2)})
		\ge u_n^{(1)}(b_n^{(1)})u_n^{(2)}(b_n^{(2)})
		>n^2>M,
		\]
		a contradiction. It proves that $u$ is not locally bounded
		above at any point of
		$\R^{N_1}\times\R^{N_2}$.
	\end{proof}

	\begin{remark}
		Restricting $u$ to any bounded product domain gives a
		counterexample to the assertion that a separately
		subharmonic function is jointly subharmonic outside
		a relatively closed nowhere dense set. In particular,
		there cannot exist relatively closed nowhere dense sets
		$S_1\subset\Omega_1$ and $S_2\subset\Omega_2$ such that
		$u$ is subharmonic on
		\[
		(\Omega_1\times\Omega_2)\setminus(S_1\times S_2).
		\]
		The function $u$ is defined globally, so the requirement
		that it be defined on a neighborhood of the closure
		of the product domain does not remove the counterexample.
	\end{remark}

	\section{Consequences for earlier results}
		In \cite{Kalantar}, we studied families of subharmonic and separately
	subharmonic functions, focusing on their local boundedness properties
	and applications to separate subharmonicity.
	It was later observed by Ahmed Zeriahi that the proof of
	\cite[Lemma~3.2]{Kalantar} contains an error. Namely, it is incorrectly
	claimed that
	\[
	v<\lambda \ \text{on }\omega
	\quad\Longrightarrow\quad
	v\geq\lambda \ \text{on }\partial\omega.
	\]
	Consequently, the decomposition
	\[
	\partial\omega
	=
	\{\zeta\in\partial\omega:v(\zeta)=\lambda\}
	\cup
	\{\zeta\in\partial\omega:v(\zeta)>\lambda\}
	\]
	and the resulting proof that
	\[
	\{\zeta\in\partial\omega:v(\zeta)>\lambda\}
	\]
	is polar are not justified. Since the proof of
	\cite[Theorem~3.1]{Kalantar} relies on this lemma, that proof is
	incomplete as well. 
	
	Worse still, the following counterexample, communicated by Vestislav
	Apostolov, shows that \cite[Lemma~3.2]{Kalantar} is false.
	\begin{counterexample}\label{counterexample}
		Let \(B\subset\R^{N}\), \(N\geq3\), be an open ball. By
		\cite[Example~3.3.2]{Armitage Gardiner book}, there exists an entire
		subharmonic function \(u\), finite almost everywhere,
		such that \(u^{-1}(\{-\infty\})\) is dense in \(B\).  Choose \(b\in B\) such that \(u(b)>-\infty\), and choose
		\(\lambda<u(b)\). Set
		\[
		\omega:=\{x\in B:u(x)<\lambda\}.
		\]
		Since \(u\) is upper semicontinuous, \(\omega\) is open. Moreover,
		\(u^{-1}(\{-\infty\})\cap B\subset\omega\), so \(\omega\) is dense in \(B\).
		Consequently,
		\[
		\operatorname{int}\overline{\omega}=B,
		\] where int designates the interior.
		Since \(b\in B\) and \(u(b)>\lambda\), this contradicts the conclusion
		of \cite[Lemma~3.2]{Kalantar}.

		Furthermore, the set
		\[
		E:=\{x\in B:u(x)>\lambda\}
		\]
		has positive Lebesgue measure. Indeed, otherwise \(u\leq\lambda\)
		almost everywhere on every ball \(B(b,r)\) such that
		\(\overline{B(b,r)}\subset B\), and the submean inequality would give
		\[
		u(b)\leq
		\frac{1}{|B(b,r)|}\int_{B(b,r)}u(x)\,dx
		\leq\lambda,
		\]
		contrary to the choice of \(b\). Since \(\omega\) is open and dense in
		\(B\), we have \(B\setminus\omega\subset\partial\omega\). Hence
		\[
		E\subset
		\{\zeta\in\partial\omega:u(\zeta)>\lambda\}.
		\]
		The latter set therefore has positive Lebesgue measure and consequently
		is not polar. Thus the conclusion of
		\cite[Lemma~3.2]{Kalantar} is false.

	\end{counterexample}

Theorem \ref{S_Can_be_R} shows that the local boundedness conclusion of \cite[Theorem 3.1]{Kalantar} fails in dimensions \(N\geq3\). Likewise, Theorem \ref{Th 7_4} shows that the corresponding conclusion of \cite[Theorem 4.1]{Kalantar} fails when both factors have dimension at least three.
	
\section*{Acknowledgments}
 I thank Ahmed Zeriahi for pointing out the error in the proof of \cite[Lemma 3.2]{Kalantar} and Vestislav Apostolov for communicating Counterexample \ref{counterexample}.

\end{document}